\documentclass[11pt]{article}
\usepackage{indentfirst}
\usepackage{amsfonts}
\usepackage{graphicx}
\usepackage{url}
\usepackage{algorithm}
\usepackage{algpseudocode}
\usepackage{listings}
\newtheorem{exemplu}{Example}[section]

\author{E. Scheiber\thanks{e-mail: scheiber@unitbv.ro}} 
\title{On the Variational Iteration Method for the Nonlinear Volterra Integral Equation}
\date{}
\begin{document}
\maketitle

\begin{abstract}
The variational iteration method is used to solve nonlinear Volterra integral equations. Two approaches
are presented distinguished by the method to compute the Lagrange multiplier. 


\end{abstract}

\section{Introduction}

The Ji-Huan He's Variational Iteration Method (VIM) was applied to a large range differential and
integral equations problems \cite{1}. The main ingredient of the VIM is the Lagrange multiplier used to improve
an approximation of the solution of the problem to be solved \cite{2}.

In this paper the VIM  to solve a nonlinear Volterra equation is resumed.
As specified in \cite{30} the Volterra integral equation must first be
transformed to an ordinary differential equation or a
nonlinear Volterra integro-differential equation by differentiating both sides.
The solution of the Volterra integral equation applying VIM was exemplified in a
series of papers \cite{34}, \cite{33}, \cite{31}, \cite{32}.

By applying the VIM, an initial value problem is deduced in order to obtain the Lagrange multiplier.
The trick is to perform such a variation that the solution of the generated initial value problem can
be deduced analytically.

We analyze two variants to compute the Lagrange multiplier.
The first variant was used in \cite{30} and \cite{31}. We observed that this approach reduces to the 
successive approximation method.

Finally two numerical examples are given showing that the second  approach gives a more rapidly converging version of VIM.

\section{The nonlinear Volterra integral equation\\ of second kind}

The nonlinear Volterra integral equation of second kind is \cite{30}
\begin{equation}\label{vim1et1}
y(x)=f(x)+\int_0^xK(x,t)F(y(t))\mathrm{d}t,\qquad x\in[0,x_f],
\end{equation}
where
\begin{itemize}
\item
$f(x), K(x,t)$ are continuous derivable functions;
\item
$F(y)$ has a continuous second order derivative.
\end{itemize}

These conditions ensure the existence and uniqueness of the solution to the nonlinear Volterra integral equation (\ref{vim1et1}).
As a consequence the solution can be computed using the successive approximation method (SAM)
\begin{eqnarray}
u_{k+1}(x)&=&f(x)+\int_0^xK(x,t)F(u_k(t))\mathrm{d}t,\label{vim1et2}\\
u_0(x)&=& f(x) \nonumber 
\end{eqnarray}
and $y(x)=\lim_{k\rightarrow\infty}u_k(x).$

\section{The variational iteration method}

The main ingredient of the VIM is the Lagrange multiplier used to improve
the computed approximations relative to a given iterative method.

\subsection*{The initial approach}
We shall remind the VIM applied to solve a nonlinear Volterra integral equation as it were presented in\cite{30}, \cite{31}.
In this approach the variation of the unknown function in the nonlinear term is neglected,
resulting in a easily solvable initial value problem for the Lagrange multiplier.

The derivative of (\ref{vim1et1}) is used in VIM for the Volterra integral equation
\begin{equation}\label{vim1et3}
u_{k+1}(x)=u_k(x)+\int_0^x\Lambda(t)\left(u'_k(t)-f'(t)-\frac{\mathrm{d}}{\mathrm{d}t}\int_0^tK(t,s)F(u_k(s))\mathrm{d}s\right)\mathrm{d}t.
\end{equation}
The variation will be not applied  to the nonlinear term $F(u_k(s))$ and the above equality is rewritten as
\begin{equation}\label{vim1et4}
u_{k+1}(x)=u_k(x)+\int_0^x\Lambda(t)\left(u'_k(t)-f'(t)-\frac{\mathrm{d}}{\mathrm{d}t}\int_0^tK(t,s)F(\tilde{u}_k(s))\mathrm{d}s\right)\mathrm{d}t.
\end{equation}
If $u_{k+1}(x)=y(x)+\delta u_{k+1}(x), u_k(x)=y(x)+\delta u_k(x),$ but $\tilde{u}_k(x)=y(x),$ then it results
\begin{equation}\label{vim1et5}
\delta u_{k+1}(x)=\delta u_k(x)+\int_0^x\Lambda(t)\delta u'_k(t)\mathrm{d}t.
\end{equation}
After an integration by parts it results
$$
\delta u_{k+1}(x)=(1+\Lambda(x))\delta u_k(x)+\int_0^x\Lambda'(t)\delta u_k(t)\mathrm{d}t.
$$
In order that $u_{k+1}$ be a better approximation than $u_k,$ it is required that $\Lambda$ is the solution
of the following initial value problem
\begin{eqnarray}
\Lambda'(x)&=&0, \label{vim1et6}\\
\Lambda(t)&=&-1,\qquad t\in[0,x].\label{vim1et7}
\end{eqnarray}
The solution of this initial value problem in $\Lambda(t)=-1.$ Replacing this solution in (\ref{vim1et3}) we get 
$$
u_{k+1}(x)=u_k(x)-\int_0^x\left(u'_k(t)-f'(t)-\frac{\mathrm{d}}{\mathrm{d}t}\int_0^tK(t,s)F(u_k(s))\mathrm{d}s\right)\mathrm{d}t=
$$
$$
=f(x)+\int_0^xK(x,s)F(u_k(s))\mathrm{d}s+\left(u_k(0)-f(0)\right)=
$$
$$
=f(x)+\int_0^xK(x,s)F(u_k(s))\mathrm{d}s.
$$
Thus we found the successive approximation method (\ref{vim1et2}).

\subsection*{Another approach}

We shall take into account partially the variation of the nonlinear term containing the unknown function.
Now in (\ref{vim1et3}) we apply  Leibniz's rule of differentiation under the integral sign 
\begin{equation}\label{vim1et12}
u_{k+1}(x)=u_k(x)+
\end{equation}
$$
+\int_0^x\Lambda(t)\left(u'_k(t)-f'(t)-K(t,t)F(u_k(t))-\int_0^t\frac{\partial K(t,s)}{\partial t}F(u_k(s))\mathrm{d}s\right)\mathrm{d}t.
$$
We require that the variation doesn't affect the term with integral, i.e.:
$$
u_{k+1}(x)=u_k(x)+
$$
$$
+\int_0^x\Lambda(t)\left(u'_k(t)-f'(t)-K(t,t)F(u_k(t))-\int_0^t\frac{\partial K(t,s)}{\partial t}F(\tilde{u}_k(s))\mathrm{d}s\right)\mathrm{d}t.
$$
It results
\begin{equation}\label{vim1et8}
\delta u_{k+1}(x)=\delta u_k(x)+
\end{equation}
$$
+\int_0^x\Lambda(t)
\left(\delta u'_k(t)-K(t,t)\left(F(u_k(t))-F(y(t))\right)\right)\mathrm{d}t=
$$
$$
=\delta u_k(x)+\int_0^x\Lambda(t)
\left(\delta u'_k(t)-K(t,t)F'(y(t))\delta u_k(t)\right)\mathrm{d}t+O((\delta u_k)^2)=
$$
$$
=(1+\Lambda(x))\delta u_k(x)-\int_0^x\left(\Lambda'(t)+\Lambda(t)K(t,t)F'(y(t))\right)\delta u_k(t)\mathrm{d}t+O((\delta u_k)^2).
$$
Again, in order for $u_{k+1}$ to be a better approximation than $u_k,$ it is required that $\Lambda$ is the solution
of the following initial value problem
\begin{eqnarray}
& \Lambda'(t)+\Lambda(t)K(t,t) F'(y(t))=0, \label{vim1et9}\\
&\Lambda(x)=-1.&\label{vim1et10}
\end{eqnarray}
for $t\in[0,x].$
The solution of this initial value problem (\ref{vim1et9})-(\ref{vim1et10}) is given by
\begin{equation}\label{vim1et11}
\Lambda(t)=-e^{\int_t^xK(s,s)F'(y(s))\mathrm{d}s}.
\end{equation}
Because $y(s)$ is an unknown function, the following problem is considered instead of (\ref{vim1et9})-(\ref{vim1et10}) 
\begin{eqnarray}
&\Lambda'(t)+\Lambda(t)K(t,t) F'(u_k(t))=0, \label{vim1et90}\\
&\Lambda(x)=-1\label{vim1et91}
\end{eqnarray}
with the solution denoted  $\Lambda_k(t).$ The solution of the problem (\ref{vim1et90})-(\ref{vim1et91}) is
$$
\Lambda_k(t)=-e^{\int_t^xK(s,s)F'(u_k(s))\mathrm{d}s}.
$$


Our next goal is to find a convenient form to implement (\ref{vim1et12}) which is equivalent with (\ref{vim1et3}).
Denoting
$$
v_{k+1}(t)=u_k(t)+\int_0^tK(t,s)F(u_k(s))\mathrm{d}t
$$
(\ref{vim1et3}) may be rewritten as
$$
u_{k+1}(x)=u_k(x)+\int_0^x\Lambda_k(t)u'_k(t)\mathrm{d}t-
$$
$$
-\int_0^x\Lambda_k(t)\left(f'(t)+\frac{\mathrm{d}}{\mathrm{d}t}\int_0^tK(t,s)F(u_k(s))\mathrm{d}s\right)\mathrm{d}t.
$$
After integration by parts in the two integrals we obtain
$$
u_{k+1}(x)=(1+\Lambda_k(x))u_k(x)-\Lambda_k(0)u_k(0)-\int_0^x\Lambda_k'(t)u_k(t)\mathrm{d}t-
$$
$$
-\Lambda_k(x)v_{k+1}(x)+\Lambda_k(0)f(0)+\int_0^x\Lambda_k'(t)v_{k+1}(t)\mathrm{d}t.
$$
Considering (\ref{vim1et10}) and that $u_k(0)=f(0)$ the above equality becomes
\begin{equation}\label{vim1et13}
u_{k+1}(x)=v_{k+1}(x)+\int_0^x\Lambda_k'(t)\left(v_{k+1}(t)-u_k(t)\right)\mathrm{d}t,
\end{equation}
where
$$
\Lambda_k'(t)=-\Lambda_k(t)K(t,t)F'(u_k(t))=e^{\int_t^xK(s,s)F'(u_k(s))\mathrm{d}s}K(t,t)F'(u_k(t)).
$$

Comparing SAM with this approach of VIM we have found experimentally that VIM needs a smaller number of iterations
to fulfill a stopping condition, e.g. the absolute error between two successive approximations is less than
a given tolerance.

\section{Numerical results}

On an equidistant mesh $x_i=i h,\ i\in\{0,1,\ldots,n\},$ with $h=\frac{x_f}{n},$ in the interval $[0,x_f],$ 
the numerical solution is $(u_i^{(k)})_{0\le i\le n},$ where
$$
u_i^{(k)}\approx u_k(x_i).
$$
Details of implementation:
\begin{itemize}
\item
The initial approximations were chosen as $u_i^{(0)}=f(x_i),\ i\in\{0,1,\ldots,n\};$
\item
$u_0^{(k)}=f(x_0),$ for any $k\in\mathbb{N};$
\item
The component $u_i^{(k+1)}$ is computed using  (\ref{vim1et13}), for $x=x_i$. The integrals were computed
using the trapezoidal rule of numerical integration;
\item
The stopping condition was $\max_{0\le i\le n}|u_i^{(k+1)}-u_i^{(k)}|<\varepsilon.$
\end{itemize}

\begin{exemplu} (\cite{30}, p. 241, Example 3)
\begin{equation}\label{vim1ex3} 
y(x)=\cos{x}-\frac{1}{4}\sin{2x}-\frac{1}{2}x+\int_0^xy^2(t)\mathrm{d}t
\end{equation}
with the solution $y(x)=\cos{x}.$
\end{exemplu}

The evolution of the computations between the SAM and the VIM are given in the next table.

\begin{tabular}{|c|c|c|c|}
\hline\hline
\multicolumn{2}{|c|}{Successive Approximation Method}&\multicolumn{2}{|c|}{Iterative Variational Method}\\
\hline
Iteration & Error & Iteration & Error  \\
\hline\hline
1 & 5.200451 & 1 & 0.986121\\
2 & 3.280284 & 2 & 0.486573 \\
3 & 1.028282 & 3 & 0.092114\\
4 & 1.167517 & 4 & 0.005001\\
5 & 0.839565 & 5 & 0.000094\\
6 & 0.493593 & 6 & 0.000001\\
\cline{3-4}
7 & 0.198224 & & \\
8 & 0.075078 & & \\
9 & 0.023424 & & \\
10 & 0.006654 & & \\
11 & 0.001705 & & \\
12 & 0.000403 & & \\
13 & 0.000089 & & \\
14 & 0.000018 & &\\
15 & 0.000004 & &\\
\hline\hline
\end{tabular}

\vspace*{0.3cm}
\noindent
In the above table, the meaning of the field \textit{Error} at iteration $k$ is given by the expression
$\max_{0\le i\le n}|u_i^{(k)}-u_i^{(k-1)}|.$

The used parameters to obtained the above results were: $x_f=\pi,\ n=30,\ \varepsilon=10^{-5}.$

\begin{figure}[htbp]
\begin{center}
\includegraphics[width=10cm,height=8cm,keepaspectratio]{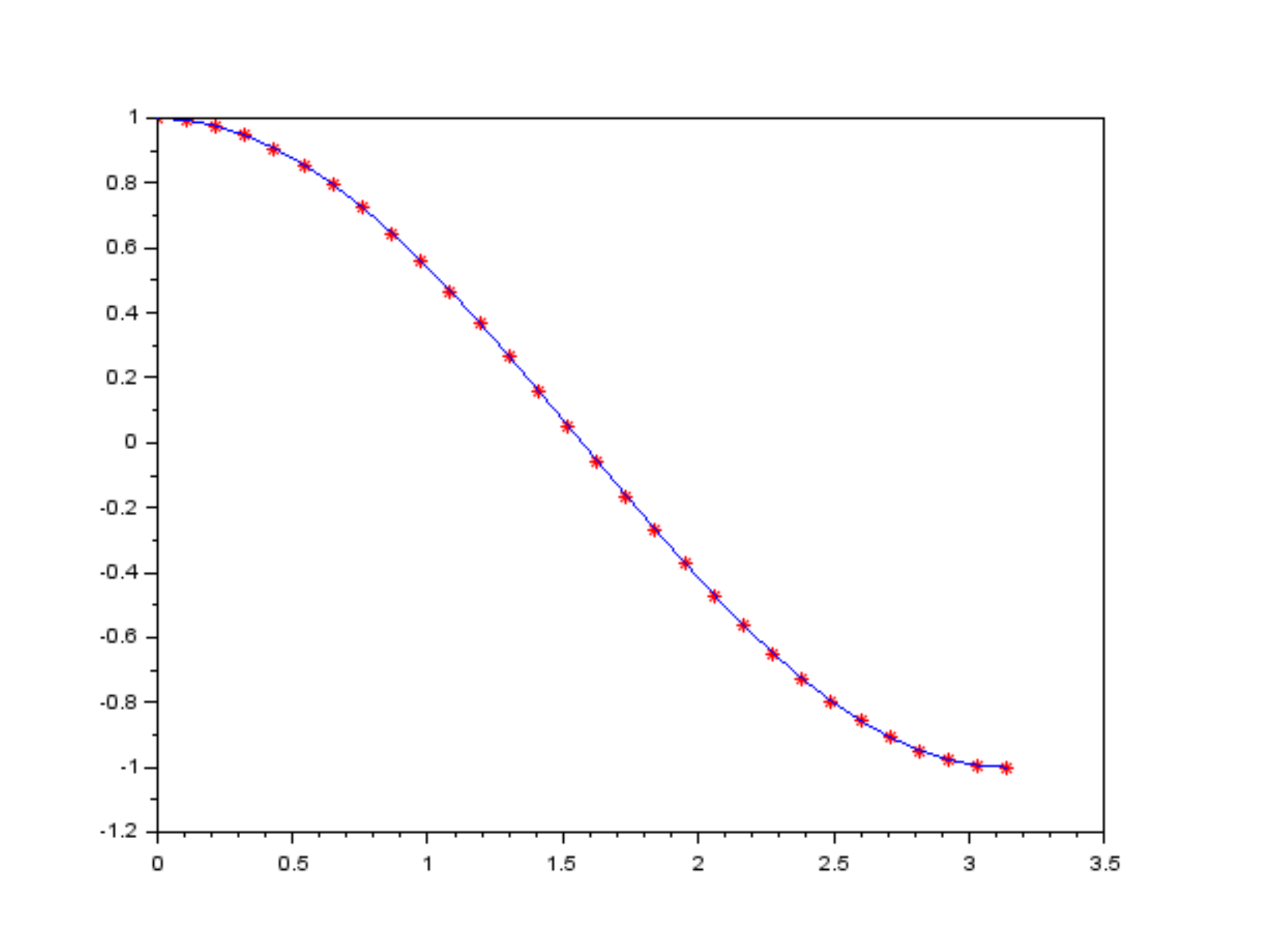}

\caption{Plot of the numerical solution and the exact solution.}\label{vim1pex3}
\end{center}
\end{figure}

The maximum of absolute errors between the obtained numerical solution and the exact solution was 0.003538.
The plot of the numerical solution vs. the solution of (\ref{vim1ex3}) are given in Fig. \ref{vim1pex3}.

\begin{exemplu} (\cite{30}, p. 240, Example 2)
\begin{equation}\label{vim1ex3} 
y(x)=e^x-\frac{1}{3}xe^{3x}+\frac{1}{3}x+\int_0^xxy^3(t)\mathrm{d}t
\end{equation}
with the solution $y(x)=e^x.$
\end{exemplu}

\begin{figure}[htbp]
\begin{center}
\includegraphics[width=10cm,height=8cm,keepaspectratio]{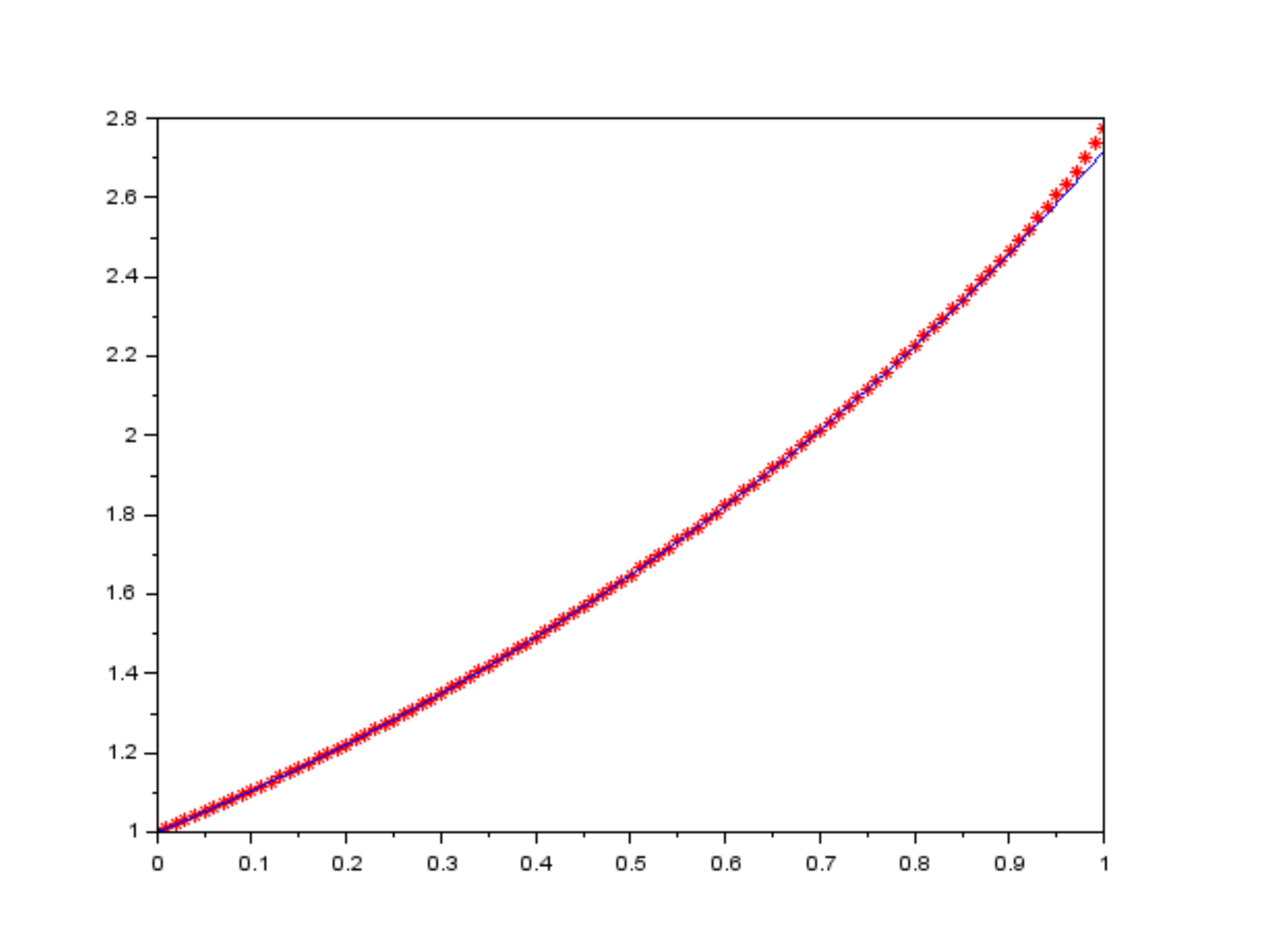}

\caption{Plot of the numerical solution and the exact solution.}\label{vim1pex2}
\end{center}
\end{figure}

The corresponding results are given in the next Table and Fig. \ref{vim1pex2}.

\begin{tabular}{|c|c|c|c|}
\hline\hline
\multicolumn{2}{|c|}{Successive Approximation Method}&\multicolumn{2}{|c|}{Iterative Variational Method}\\
\hline
Iteration & Error & Iteration & Error  \\
\hline\hline
1 & 1.314327 & 1 & 4.085689\\
2 & 0.990523 & 2 & 1.844170 \\
3 & 1.791381 & 3 & 1.202564\\
4 & 1.733644 & 4 & 1.073406\\
5 & 0.732878 & 5 & 1.111831\\
6 & 0.464864 & 6 & 0.700407\\
7 & 0.424189 & 7 & 0.13054\\
8 & 0.409519 & 8  & 0.004592 \\
9 & 0.398283 & 9 & 0.00009\\
10 & 0.381286 & 10 & 0.000002\\
\cline{3-4}
11 & 0.352116 & & \\
12 & 0.307541 & & \\
\vdots & \vdots & & \\
25 & 0.000011 & &\\
26 & 0.000003 & &\\
\hline\hline
\end{tabular}

\vspace*{0.3cm}
\noindent
The used parameters to obtained the results were: $x_f=1.0,\ n=101,\ \varepsilon=10^{-5}.$

The maximum of absolute errors between the obtained numerical solution and the exact solution was 0.056763.

\section{Conclusions}
\begin{enumerate}
\item
The VIM as it is presented in \cite{30}, \cite{31}, when $\Lambda(x)=-1,$ for any $x,$ 
returns to the successive approximation method.
\item
Comparing SAM with VIM we have found experimentally that VIM needs a smaller number of iterations
to fulfill a stopping condition, e.g. the absolute error between two successive approximations is less than
a given tolerance. 
\end{enumerate}

\end{document}